\documentclass[11pt]{article}
\usepackage{geometry}                
\usepackage{pgf,tikz}
\usetikzlibrary{automata,arrows}
\usepackage{graphicx}
\usepackage{amssymb}
\usepackage{epstopdf}
\title{Permutations sortable by two stacks in series}
\author{Andrew Elvey-Price and Anthony J Guttmann}
\date{}                                           

\begin{document}
\maketitle
\begin{abstract}
We address the problem of the number of permutations that can be sorted by two stacks in series. We do this by first counting all such permutations of length less than 20 exactly, then using a numerical technique to obtain nineteen further coefficients approximately. Analysing these coefficients by a variety of methods we conclude that the OGF behaves as
$$S(z) \sim A (1 - \mu \cdot z)^\gamma,$$ where $\mu =12.45 \pm 0.15,$  $\gamma= 1.5 \pm 0.3,$ and $A \approx 0.02$.

\end{abstract}
\section{Introduction}
In the late 1960s Knuth \cite{K68} introduced the idea of classifying the common data structures of computer science in terms of the number of permutations of length $n$ that could be sorted by the given data structure, to produce the identity permutation. Knuth demonstrated the usefulness of this approach by showing that a simple stack could sort all such permutations except those which had any three elements in relative order 312. This restriction meant that of the $n!$ possible permutations of length $n,$ only $C_n \sim 4^n/(n^{3/2}\sqrt{\pi})$ could be sorted by a simple stack. Here $C_n$ denotes the cardinality of the $n$th Catalan number. Knuth went on to pose the same question for more complex data structures, such as a double-ended queue or {\em deque}, which is a linear list in which insertions and deletions can take place at either end. In a later volume of his celebrated book \cite{K73}, he asked the same question about compositions of stacks.

\begin{figure}[h!] 
   \centering
   \includegraphics[width=4in]{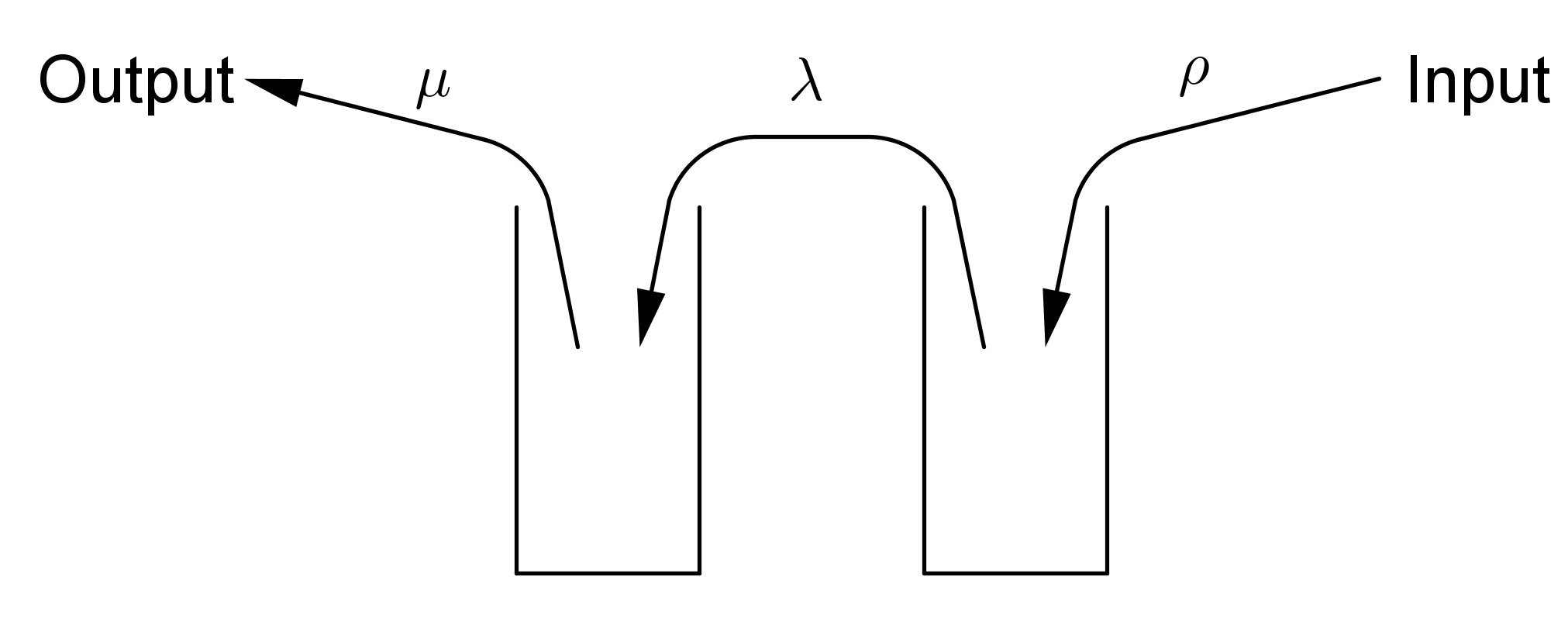} 
   \caption{Two stacks in series}
   \label{fig:tsis}
\end{figure}

The three most interesting, and most intensively studied permutation-related sorting problems associated with data structures relate to permutations that can be sorted by (i) a deque, (ii)  two stacks in parallel (2SIP) and (iii) two stacks in series (2SIS). The data structure corresponding to two stacks in series is shown in Fig. \ref{fig:tsis}. A permutation of length $n$ is said to be {\em sortable} if it is possible to start with this permutation as the input, and output the numbers $1,2,\ldots, n$ in order, using only the moves $\rho$, $\lambda$ and $\mu$ in some order. Here $\rho$ pushes the next element from the input onto the first stack, $\lambda$ pushes the top element of the first stack onto the top of the second stack, and $\mu$ outputs (pops) the top element of the second stack to the output stream, as shown in Fig. \ref{fig:tsis}.

Recently Albert and Bousquet-M\'elou \cite{BMA14} solved the problem relating to two stacks in parallel, while subsequently we \cite{EPG15} related the solution of the 2SIP problem to the solution of the deque problem. This leaves only the 2SIS problem unresolved. Significant progress has been made on subsets of that problem. For example Atkinson, Murphy and Ru{\v{s}}kic \cite{AMR02} solved the problem in the case of {\em sorted} stacks, while Elder, Lee and Rechnitzer  \cite{ELR15} solved the problem in the case when one of the stacks is of depth 2. Unfortunately, both these cases correspond to an exponentially small subset of the full set of stack-sortable permutations. In \cite{PR13}, Pierrot and Rossin give a polynomial algorithm to decide if a given permutation is sortable by two stacks in series.

In all cases we've mentioned, the number of permutations of length $n$ that can be sorted by the given data structure grows exponentially (just as in the simple stack case discussed above), and more precisely, it is expected that $p_n,$ the number of permutations of length $n$ sortable by any of the afore-mentioned data structures, behaves as $p_n \sim const \cdot  \mu^n \cdot n^g$ in general.  In \cite{AAL10}, rigorous upper and lower bounds on $\mu$ are given for deque sorting, and also for 2SIP and 2SIS. For 2SIS the bounds are $8.156 < \mu < 13.374.$

In this paper we give an alternative approximation. We have evaluated the exact number of stack-sortable permutations of length $n$ for $n < 20,$ and describe numerical techniques that give the approximate number for $20 \le n \le 30.$ We then apply a range of standard and specialised series analysis techniques \cite{G89} to conjecture the asymptotics of the generating function coefficients. This is a computationally difficult problem. The only existing series we can find are for $n \le 10$ in the PhD thesis of Pierrot \cite{P13}, and the last two values are incorrect.

If $$S(z) = \sum_{n \ge 0} s_n \cdot z^n$$ is the ordinary generating function for the number of permutations sortable by two stacks in series, then we find
$$S(z) \sim A (1 - \mu \cdot z)^\gamma,$$ where $\mu \approx 12.4,$  $\gamma \approx 1.5,$ and $A \approx 0.026.$

In the next section we describe the derivation of the coefficients $s_n,$ and in the subsequent section we give our analysis of the data.

\section{Generating coefficients of OGF}

\subsection{Basic algorithm}

We start with a simple, but inefficient algorithm to calculate the coefficients of the OGF, on which our more efficient algorithm is based. Consider the three moves $\rho$, which pushes the next element from the input onto the first stack, $\lambda$, which pushes the top element of the first stack onto the second stack, and $\mu$, which outputs the top element of the second stack as shown in Fig. \ref{fig:tsis}. We have already defined sortable permutations. We call a permutation of length $n$ {\em achievable} if it is possible to output that permutation, starting with the numbers $1,2,\ldots, n$ in order.  Rather than enumerating sortable permutations directly, we will instead enumerate achievable permutations, since the two classes share the same OGF. We call a word $w$ over the alphabet $\{\rho,\lambda,\mu\}$ an {\em operation sequence} if $w$ corresponds to a permutation.  That is, $w$ is called an operation sequence if $w$ contains an equal number of occurrences of each of the three letters, and after any point in $w$, the letter $\rho$ has appeared at least as many times as $\lambda$, which has appeared at least as many times as $\mu$. Call two operation sequences equivalent if they produce the same permutation. Note that this also means that they sort the same permutation. The basic algorithm, which we will call algorithm 1, works as follows:
\begin{itemize}
\item Define the function $addreachableperms$ which takes in the state $S$ of the sorting machine, and a set of permutations and adds every permutation which can be achieved from that state to the set, by recursively calling the same function on each of the three or fewer states which can be reached from $S$ by one of the moves $\rho$, $\mu$, or $\gamma$.
\item create an empty set $P$ of permutations.
\item Call the function $addreachableperms$ on the initial state of the stack and the set $P$
\item Then the $n$th coefficient of the OGF is equal to the size of $P$, since the permutations in $P$ are exactly the achievable permutations of size $n$.
\end{itemize}

This algorithm is very slow because it has to consider all operation sequences of size $3n$ separately, and the number of operation sequences of length $3n$ grows like $27^n$.

\subsection{Forbidden words and regular languages}
The first improvement which we make is to reduce the number of operation sequences which the algorithm has to consider by removing many operation sequences which create the same permutation. Call two operation sequences equivalent if they create the same permutation. We define the ordering $\rho<\lambda<\mu$, and we call an operation sequence optimal if it is lexicographically larger than any other equivalent operation sequence. Rather than parse all operation sequences of size $3n$, we now only insist that we parse all optimal operation sequences, since these will still create all achievable permutations. Call a word $v$ over the alphabet $\{\rho,\lambda,\mu\}$ forbidden if there is another word $v'>v$, which has the same effect on the sorting machine as $v$. Note that if an operation sequence $w$ has a forbidden subword $v$, then we can change $v$ to $v'$ in $w$ to create an equivalent operation sequence $w'$. Moreover, $w'>w$, so $w$ is not optimal. Hence, any optimal operation sequence contains no forbidden words. Note that $\rho\mu$ is a forbidden word, since it has the same effect on the sorting machine as $\mu\rho$. Also, $\rho\lambda\mu\lambda$ is a forbidden word since it has the same effect on the sorting machine as $\lambda\rho\lambda\mu$. For letters $x$ and $y$, we call a word $v$ over the alphabet $\{x,y\}$ an $x,y$-Catalan word if the following conditions hold:
\begin{itemize}
\item $v$ contains an equal number of $x$'s and $y$'s
\item for any leading subword $u$ of $v$, the word $u$ contains at least as many $x$'s as $y$'s.
\end{itemize}
In other words, if we replace each $x$ in $v$ with an up step and each $y$ in $v$ with a down step, we get a Dyck path. Note that if $u$ is a $\rho,\lambda$-Catalan word, and $v$ is a $\lambda,\mu$-Catalan word, then the effect of $u$ on the sorting machine is to move and permute items from the input to the second stack. The effect of $v$ is to move and permute items from the first stack to the output. Hence, these two operations commute, so $uv$ and $vu$ are equivalent. Since $u$ begins with $\rho$ and $v$ begins with $\lambda$, we have $uv<vu$, so $uv$ is a forbidden word.

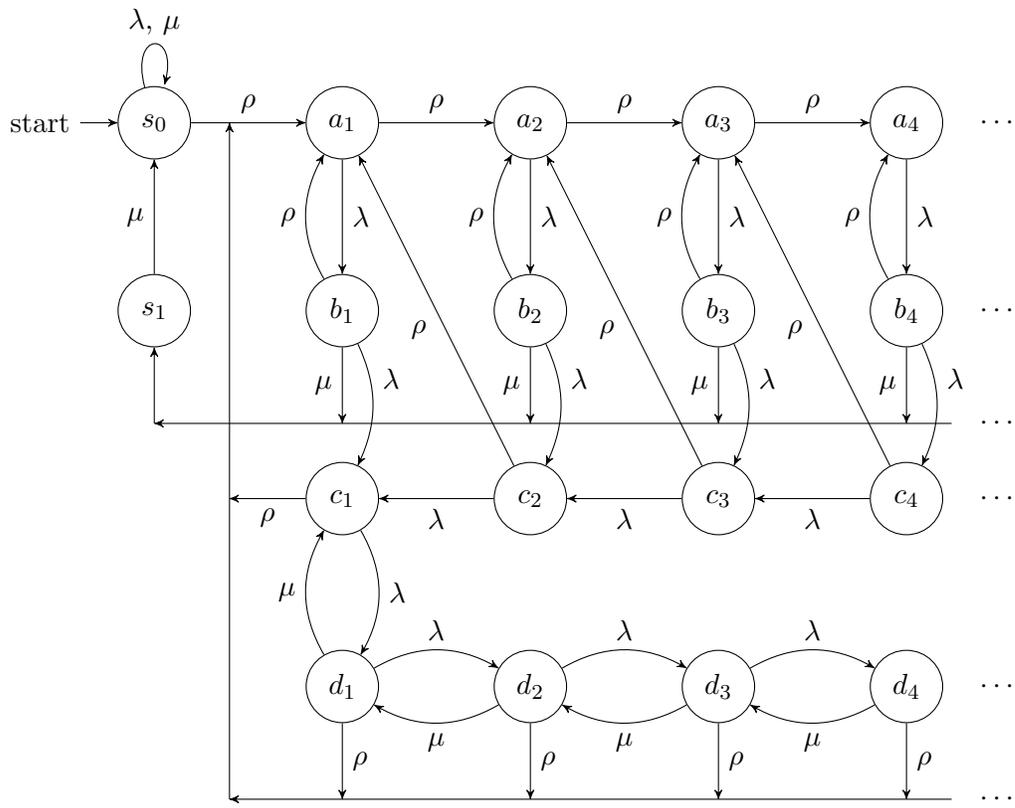
\begin{figure}
\begin{tikzpicture}[>=stealth',shorten >=0pt,auto,node distance=2.5 cm, scale = 1, transform shape]

\node[initial,state] (X)                                    {$s_0$};
\node[state]         (Y) [below of=X]                       {$s_1$};
\node[state]         (A1) [right of=X]                       {$a_1$};
\node[state]         (B1) [below of=A1]                       {$b_1$};
\node[state]         (C1) [below of=B1]                       {$c_1$};
\node[state]         (D1) [below of=C1]                       {$d_1$};
\node[state]         (A2) [right of=A1]                       {$a_2$};
\node[state]         (B2) [below of=A2]                       {$b_2$};
\node[state]         (C2) [below of=B2]                       {$c_2$};
\node[state]         (D2) [below of=C2]                       {$d_2$};
\node[state]         (A3) [right of=A2]                       {$a_3$};
\node[state]         (B3) [below of=A3]                       {$b_3$};
\node[state]         (C3) [below of=B3]                       {$c_3$};
\node[state]         (D3) [below of=C3]                       {$d_3$};
\node[state]         (A4) [right of=A3]                       {$a_4$};
\node[state]         (B4) [below of=A4]                       {$b_4$};
\node[state]         (C4) [below of=B4]                       {$c_4$};
\node[state]         (D4) [below of=C4]                       {$d_4$};
\node[draw=none]   at (11.2,0)                       {$\cdots$};
\node[draw=none]   at (11.2,-2.5)                       {$\cdots$};
\node[draw=none]   at (11.2,-5)                       {$\cdots$};
\node[draw=none]   at (11.2,-7.5)                       {$\cdots$};
\node[draw=none]   at (11.2,-4)                       {$\cdots$};
\node[draw=none]   at (11.2,-9)                       {$\cdots$};

\path[->] (X) edge   [loop above]     node [align=center]  {$\lambda$, $\mu$} (X)
      (X) edge      node  {$\rho$} (A1)
      (Y) edge      node [align=center]  {$\mu$} (X)
      (C1) edge  [bend left]    node [align=center]  {$\lambda$} (D1)
      (D1) edge  [bend left]    node [align=center]  {$\mu$} (C1)

      (A1) edge       node [align=center]  {$\rho$} (A2)
      (A1) edge     node [align=center]  {$\lambda$} (B1)
      (B1) edge   [bend left]   node [align=center]  {$\rho$} (A1)
      (B1) edge   [bend left=25]   node [pos=0.45] {$\lambda$} (C1)
      (C2) edge      node [align=center]  {$\lambda$} (C1)
      (C2) edge   node   {$\rho$} (A1)
      (D1) edge  [bend left]    node [align=center]  {$\lambda$} (D2)
      (D2) edge  [bend left]    node [align=center]  {$\mu$} (D1)

      (A2) edge       node [align=center]  {$\rho$} (A3)
      (A2) edge     node [align=center]  {$\lambda$} (B2)
      (B2) edge   [bend left]   node [align=center]  {$\rho$} (A2)
      (B2) edge   [bend left=25]   node [pos=0.45] {$\lambda$} (C2)
      (C3) edge      node [align=center]  {$\lambda$} (C2)
      (C3) edge      node [align=center]  {$\rho$} (A2)
      (D2) edge  [bend left]    node [align=center]  {$\lambda$} (D3)
      (D3) edge  [bend left]    node [align=center]  {$\mu$} (D2)

      (A3) edge       node [align=center]  {$\rho$} (A4)
      (A3) edge     node [align=center]  {$\lambda$} (B3)
      (B3) edge   [bend left]   node [align=center]  {$\rho$} (A3)
      (B3) edge   [bend left=25]   node [pos=0.45] {$\lambda$} (C3)
      (C4) edge      node [align=center]  {$\lambda$} (C3)
      (C4) edge      node [align=center]  {$\rho$} (A3)
      (D3) edge  [bend left]    node [align=center]  {$\lambda$} (D4)
      (D4) edge  [bend left]    node [align=center]  {$\mu$} (D3)

      (A4) edge     node [align=center]  {$\lambda$} (B4)
      (B4) edge   [bend left]   node [align=center]  {$\rho$} (A4)
      (B4) edge   [bend left=25]   node [pos=0.45] {$\lambda$} (C4)

(C1) edge node [align=center]  {$\rho$} (1,-5)
(D2) edge node [align=center]  {$\rho$} (5,-9)
(D3) edge node [align=center]  {$\rho$} (7.5,-9)
(D4) edge node [align=center]  {$\rho$} (10,-9)
(D1) edge node [align=center]  {$\rho$} (2.5,-9)
(10.6,-9) edge (1,-9)
(1,-9) edge (1,0)

(B1) edge node [left]  {$\mu$} (2.5,-4)
(B2) edge node [left]  {$\mu$} (5,-4)
(B3) edge node [left]  {$\mu$} (7.5,-4)
(B4) edge node [left]  {$\mu$} (10,-4)
(10.6,-4) edge (0,-4)
(0,-4) edge (Y);

\path[-]
;

\end{tikzpicture}
\caption{The infinite state automaton $\Gamma$}
\label{fig:automaton}
\end{figure}

 We now construct the deterministic infinite state automaton $\Gamma$ shown in Fig. \ref{fig:automaton}, which accepts all words which are not forbidden. Note that $\Gamma$ also accepts some words which are forbidden.
For an operation sequence $w$ of size at most $3n$, the word $w$ is accepted by $\Gamma$ if and only if $w$ does not contain any of the words $\rho\mu$, $\rho\lambda\mu\lambda$ or any word of the form $uv$, where $u$ is a $\rho,\lambda$-Catalan word, and $v$ is a $\lambda,\mu$-Catalan word. Since all of these words are forbidden, any operation sequence $w$ which is not forbidden is accepted by $\Gamma$. For any integer $m$, at least $m$ occurrences of the letter $\rho$ are required to reach any of the states $a_{m}$, $b_{m}$ or $c_{m}$, and at least $m+2$ occurrences of the letter $\lambda$ are required to reach the state $d_{m}$. hence, for operation sequences of size $n$, we only need to construct the finite state automaton $\Gamma_{n}$, consisting of the $4n$ states $s_{0},s_{1},a_{1},\ldots,a_{n},b_{1},\ldots,b_{n},c_{1},\ldots,c_{n},d_{1},\ldots,d_{n-2}$.

Note that all forbidden words of length less than 7 are rejected by $\Gamma$. For each forbidden word $u$ of length at most 9, which contains no other forbidden words, we construct the DFA $\Gamma_{u}$ which accepts the language of all words which do not contain $u$. Then every word which is not accepted by $\Gamma_{u}$ contains $u$ and is hence forbidden. We then take each of these DFA's $\Gamma_{u}$ along with $\Gamma_{n}$, and construct the DFA $\Gamma_{n}'$, which accepts the intersection of the languages accepted by all of the other DFA's.

Now, our new algorithm works  as before, except that the function $addreachableperms$ also takes in the current state $A$ of $\Gamma_{n}'$, and only recursively calls itself using one of the letters which is accepted from state $A$. Now, rather than considering all operation sequences of size $3n$, the new algorithm only considers those operation sequences which are accepted by $\Gamma_{n}'$.

The improvements to the algorithm so far significantly decrease the exponential factor in the time requirement, from 27 to about 13. However the algorithm still stores every achievable permutation of length $n$ in memory at the same time. In the next section we see that the number of such permutations is approximately $12.4^n$, so any improvements of the form which we have presented so far will not reduce the time or memory requirements below this factor.

\subsection{Increment avoiding permutations}

Our next improvement to the algorithm decreases the exponential factor by 1, and we do not improve on the factor for time any more than this. Let $p=a_{1}\ldots a_{n}$ be a permutation and let $I\subset\{1,2,\ldots,n\}$. We define the subpermutation $p|_{I}$ to be the pattern of the elements from $I$ in $p$, and define $\hat{p}|_{I}=p|_{\{1,\ldots,n\}\setminus I}$. For example, $\hat{24315}|_{\{2,3\}}=24315|_{\{1,4,5\}}=213$. Note that any subpermutation of an achievable permutation is also achievable. Let $p$ be a permutation such that $a_{j+1}=a_{j}+1$ for some $j$. Since $p$ is achievable, $\hat{p}|_{\{a_{j}\}}$ is achievable. On the other hand, if $\hat{p}|_{\{a_{j}\}}$ is achievable, then there is some operation sequence $w$ which creates it. Now replace each letter which moves $a_{j}$ in $w$ with two copies of that letter, to form a new word $w'$. Then the two copies of each letter will move $a_{j}$ and $a_{j}+1$, and $a_{j}$ will enter the first stack immediately before $a_{j}+1$, then $a_{j}+1$ will enter the second stack immediately before $a_{j}$ and finally, $a_{j}$ is output immediately before $a_{j}+1$. Since the order of everything else stays the same, the word $w'$ creates the permutation $p$. Therefore, $p$ is achievable if and only if $\hat{p}|_{\{a_{j}\}}$ is achievable.

Now, instead of considering all achievable permutations with the algorithm, we only consider permutations $a_{1}\ldots a_{n}$ for which there is no $j$ such that $a_{j+1}=a_{j}+1$. Call these increment avoiding permutations. Let $t_{n}$ be the number of these permutations, and define the generating function $T(x)=t_{0}+xt_{1}+\ldots$. Then we can uniquely create any achievable permutation by choosing a permutation $q$ counted by $T$ and replacing each number in $q$ with any positive number of consecutive integers. Hence $S(x)=T(x/(1-x))$. By taking the coefficient for $x^{n}$ on both sides of this equation, we deduce that
\[s_{n}=\sum_{i=1}^{n}{n-1\choose i-1}t_{i}.\]
The only change we make to the algorithm presented previously to instead calculate the number of increment avoiding permutations, is to forbid an item from being output if it is exactly one greater than the previous item output.

\subsection{Memory consumption and parallelisation}

Using the algorithm described so far, it is still necessary to list every achievable, increment avoiding permutation of size $n$ at the same time. To avoid this restriction, we choose some positive integer $m<n$, and write a function $numpermswithstartsequence$, which inputs $n$ and a sequence $s$ of $m$ distinct elements of $\{1,\ldots,n\}$ and outputs the number of achievable, increment avoiding permutations of size $n$ which begin with the sequence $s$. This algorithm works in the same way as before except that the first $m$ elements output must be the correct elements of $s$. We then run this function on all such sequences $s$ and add up the results. For small values of $m$ this algorithm only takes a little longer than the original algorithm because most of the time is spent while the operation sequence is long and the output is nearly complete. Since we call the function on different sequences $s$ separately, it is only necessary to store all of the (achievable, increment avoiding) permutations which begin with some sequence $s$ at any one time. Note also that we only have to remember the last $n-m$ elements of each permutation. As a result, the limiting factor for this algorithm is now the time requirement.

We now parrallelise the algorithm, by running the $numpermswithstartsequence$ on different sequences $s$ at the same time on different cores. 

\subsection{Results}

We ran this algorithm for $n<20$ using $m=6$. The program ran for 43 days on 64 cores. The coefficients of the OGF for $n < 20$ are given as a list below.\\

\noindent
[1, 1, 2, 6, 24, 120, 720, 5018, 39374, 337816, 3092691, 29659731, 294107811, 2988678546, \\
30935695794, 324832481490, 3450158410649, 36993206191004, 399827092167771, \\
4351269802153188].\\

\section{Series Analysis}
\subsection{Series extension and subsequent analysis.}
We have obtained approximate values of the next nineteen coefficients, effectively doubling the length of the series,
which are sufficiently accurate to be used in the ratio analysis we describe below. Our method for obtaining these approximate values uses {\em differential approximants} \cite{G89}, which are linear, inhomogeneous ODE's of 2nd, 3rd or 4th order, constructed to yield all the exactly known coefficients in the series expansion under consideration. By varying the degrees of the polynomials multiplying each derivative, as well as the degree of the inhomogeneous polynomial, we can construct a family of such approximants. Because every differential approximant (DA) that uses all the available series coefficients implicitly predicts
all subsequent coefficients, we can calculate, approximately, all subsequent coefficients. Of course the accuracy of these predicted coefficients decreases as the order of the predicted coefficients increases, but, as we show by example below, we can get useful estimates of the next  nineteen or so coefficients. 

For every DA using all known
coefficients, we generated the subsequent nineteen coefficients. We take the mean of the predicted coefficients, with the outlying 10\% or 15\% of estimates rejected, as our estimate. We quote one standard deviation as the error. That is to say, assume we know the coefficients $a_n$ for $n \in [0,N_{max}].$ We then predict the coefficients $a_{N_{max}+1},\,a_{N_{max}+2},\,\cdots,a_{N_{max}+19}.$ Our estimate of each such coefficient is given by the mean of the values predicted by the differential approximants. We reject obvious outliers, by discarding the top and bottom 10\% of estimates. Not surprisingly, we find the smallest error is predicted for $a_{N_{max}+1},$ with the error slowly increasing as we generate further coefficients.

These predicted coefficients are well-suited to ratio type analyses, as discrepancies in say the seventh or eighth significant digit will not affect the ratio analysis in the slightest. This is particularly useful in those situations where we suspect there might be a turning point in the behaviour of ratios or their extrapolants with our exact coefficients, as these approximate coefficients are more than accurate enough to reveal such behaviour, if it is present.

As an indication of the validity of this method, we give two applications. In the first, we take the series for two stacks in parallel, for which we actually have more than 1000 coefficients \cite{EPG15}, but assume we only have the first 20 coefficients, just as in the present case for the generating function of two stacks in series. In Table \ref{tab:par-pred} we show a selection of the estimated coefficients $p(20)$ to $p(38).$ It can be seen that we predict the next coefficient with an accuracy of 13 digits, decreasing to 7 digit accuracy for the last predicted coefficient. In every case the actual error is seen to be less than one standard deviation, indeed, it is typically $1/3$ of a standard deviation.

\begin{table}[htbp]
   \begin{center}
\begin{tabular}{|p{30pt} p{140pt} p{70pt} p{70pt}|}
\hline
$N$ &$\,\,\,\,\,\,\,\,\,\,    p_N$ estimate &1 std. devn. & actual error  \\
\hline
&& &\\
20 & $1.36000505625858 \times 10^{14} $& 81 &28\\
21 &	$9.90406677134907 \times 10^{14} $& 5285 & 1778 \\
22& $7.258100272044 \times 10^{15} $& 187074 &64212\\
23 &	$5.349517582877 \times 10^{16}$ & 4807109 &1358815\\
24 &	$3.9634005851\times 10^{17} $ &$ 9.9546 \times 10^{7}$ & $3.924 \times 10^{7}$\\
25 &	$2.95046460646 \times 10^{18} $ & $1.7832\times 10^{9}$ & $5.767 \times 10^8$\\
29  & $9.435573118 \times 10^{21} $ &$7.2824 \times 10^{13}$& $2.222 \times 10^{13}$\\
32 & $4.15469546597  \times 10^{24} $ &$1.1960 \times 10^{17}$& $3.569 \times 10^{15}$\\
35 & $1.873198683303  \times 10^{27} $ &$1.5136 \times 10^{20}$& $4.030 \times 10^{19}$\\
38 & $ 8.613038855  \times 10^{29} $ &$7.2824 \times 10^{22}$&$2.222 \times 10^{21}$\\
\hline
  \end{tabular}
   \caption{Series coefficients $p_N$ for two stacks in parallel. Approximate, predicted coefficients $p_{20}$ to $p_{38}$ from several 4th order inhomogeneous DAs, the estimated and exact error.  }
   \label{tab:par-pred}
\end{center}
\end{table}

As a second demonstration of this method, assume we only have 19 terms in the generating function for two stacks in series, and we'll predict the next  coefficient, which is the last known coefficient. 
The predictions produced by fourth-order DAs are  averaged, deleting the top and bottom 10\% of estimates. In this way we estimate $s_{19} = 4.351269803411739 \times 10^{15}.$
The correct answer is $4351269802153188,$ which is estimated with an error of 1 part in the 10th significant digit by the differential approximants. The standard deviation of the estimates is 7979922,
which is six times the actual error.

In an identical manner to that described above to estimate the coefficients of the two stacks in parallel series, of course using the exact value of $s_{19},$ we have obtained estimates of the next 19 coefficients. These are given in table \ref{tab:ser-pred} below. We also give the standard deviation of the estimates, and based on the examples already discussed, we expect coefficient errors to be less than this. It can be seen that fewer significant digits are predicted than for the two stackc in parallel  series--typically 3 or 4 fewer digits at each order. Nevertheless, the precision (4 significant digits at worst, is sufficient for a simple ratio plot.

\begin{table}[htbp]
   \begin{center}
\begin{tabular}{|p{30pt} p{140pt} p{70pt}|}
\hline
$n$ &$\,\,\,\,\,\,\,\,\,\,    s_n$ estimate &1 std. devn.   \\
\hline
&& \\
20 & $4.764211695346 \times 10^{16} $& $9.207  \times 10^{6}$ \\
21 &	$5.24460896431 \times 10^{17} $&$ 9.83  \times 10^{8} $ \\
22& $5.8016808762 \times 10^{18} $& $7.962  \times 10^{10} $\\
23 &	$6.446525027 \times 10^{19}$ &$ 2.241  \times 10^{12}$ \\
24 &	$7.192361922\times 10^{20} $ &$ 7.34 \times 10^{13}$ \\
25 &	$8.05485154 \times 10^{21} $ & $2.05\times 10^{15}$ \\
26  & $9.05248613 \times 10^{22} $ &$5.10 \times 10^{16}$\\
27  & $1.02070684 \times 10^{24} $ &$1.17 \times 10^{18}$\\
28  & $1.15442858 \times 10^{25} $ &$2.47 \times 10^{19}$\\
29  & $1.30944006 \times 10^{26} $ &$4.95 \times 10^{20}$\\
30  & $1.4893068 \times 10^{27} $ &$9.38 \times 10^{21}$\\
31  & $1.6982322 \times 10^{28} $ &$1.17 \times 10^{23}$\\
32  & $1.941173 \times 10^{29} $ &$2.98 \times 10^{24}$\\
33 & $2.2239807  \times 10^{30} $ &$5.05 \times 10^{25}$\\
34 & $2.5535645  \times 10^{31} $ &$8.36 \times 10^{26}$\\
35 & $2.938088 \times 10^{32} $ &$1.345 \times 10^{28}$\\
36 & $3.387209  \times 10^{33} $ &$2.106 \times 10^{29}$\\
37 & $3.91235  \times 10^{34} $ &$3.219 \times 10^{30}$\\
38 & $ 4.52711 \times 10^{35} $ &$4.945 \times 10^{31}$\\
\hline
  \end{tabular}
   \caption{Series coefficients $s_n$ for two stacks in series. Approximate, predicted coefficients $p_{20}$ to $p_{38}$ from several 4th order inhomogeneous DAs and the estimated error.  }
   \label{tab:ser-pred}
\end{center}
\end{table}

If we wish to plot the ratios, we can do better by extrapolating the sequence of ratios produced from the coefficients predicted. That is to say, for each approximating differential approximant one calculates the ratio of successive coefficients and averages these across all differential approximants using all known coefficients - as usual discarding the uotlying 10\% or 15\% of entries. As shown in \cite{G16} this generally gives more accurate ratios than taking ratios of predicted coefficients. In this way we have obtained the next 30 ratios, and these are shown in Table \ref{tab:rat-pred}. We see that we have 10 digit accuracy in the first predicted ratio, decreasing to 4 digit accuracy in the 30th predicted ratio.

\begin{table}[htbp]
   \begin{center}
\begin{tabular}{|p{30pt} p{140pt} p{90pt}|}
\hline
$n$ &$\,\,\,\,\,\,\,\,\,\,    r_n$ estimate &1 std. devn.   \\
\hline
&& \\
20&1.094901468298879e+01& 2.11772356e-09\\
21&1.100834576534045e+01& 1.85143285e-08\\
22&1.106218007359570e+01 &8.68930350e-08\\
23&1.111147816281692e+01 &2.96091135e-07\\
24&1.115695963791318e+01& 8.01235631e-07\\
25&1.119917449576340e+01 &1.85222167e-06\\
26&1.123855154825218e+01 &3.82169883e-06\\
27&1.127543165860235e+01& 7.19040042e-06\\
28&1.131009082195476e+01& 1.25541540e-05\\
29&1.134275634825540e+01 &2.05420282e-05\\
30&1.137361924070142e+01 &3.24451987e-05\\
31&1.140284018941137e+01 &4.79461475e-05\\
32&1.143055939166936e+01 &6.84547697e-05\\
33&1.145689689567404e+01 &9.36110721e-05\\
34&1.148196167309631e+01 &1.27067931e-04\\
35&1.150584576102606e+01 &1.65506944e-04\\
36&1.152863546691212e+01 &2.12677243e-04\\
37&1.155040666870830e+01 &2.66928864e-04\\
38&1.157122649222079e+01 &3.29184438e-04\\
39&1.159116101123194e+01 &3.99718853e-04\\
40&1.161026217312874e+01 &4.79900592e-04\\
41&1.162858140362701e+01 &5.64407415e-04\\
42&1.164616370369958e+01 &6.53792233e-04\\
43&1.166306461316009e+01 &7.55633579e-04\\
44&1.167930152412725e+01 &8.55570943e-04\\
45&1.169494749603971e+01 &9.75376673e-04\\
46&1.171001288824953e+01 &1.09733642e-03\\
47&1.172454051929511e+01 &1.22834662e-03\\
48&1.173852889188502e+01 &1.35210096e-03\\
49 &1.175203710982201e+01 &1.48442613e-03\\

\hline
  \end{tabular}
   \caption{Predicted ratios $r_n=s_n/s_{n-1}$ and their standard deviation for two stacks in series, from 4th order inhomogeneous DAs.  }
   \label{tab:rat-pred}
\end{center}
\end{table}


\subsection{Series analysis of extended series}

\begin{figure}[h!] 
   \centering
   \includegraphics[width=4in]{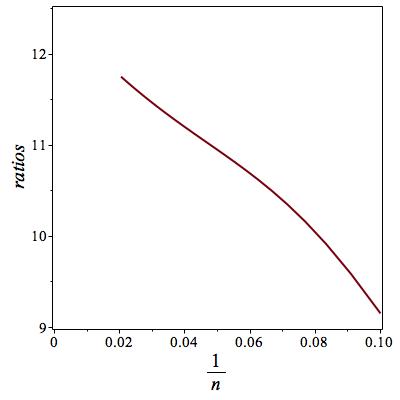} 
   \caption{Plot of ratios against $1/n.$}
   \label{fig:r1}
\end{figure}

We first performed a simple ratio analysis, under the assumption that the coefficients behave as $s_n \sim const \cdot \mu^n \cdot n^g.$ Then the ratio of successive coefficients, $r_n$  behaves as $$r_n = \frac{s_n}{s_{n-1}} = \mu \left ( 1 + \frac{g}{n} + o\left ( \frac{1}{n} \right ) \right ),$$ so plotting the ratios against $1/n$ should, for sufficiently large $n,$ give a straight line intercepting the abscissa at $\mu,$ and with gradient $g \cdot \mu.$ We show this plot in Fig. \ref{fig:r1}. One sees some low $n$ curvature, and this suggests the presence of a confluent singularity. That is to say, the generating function probably behaves as
$$S(z) \sim A(1-\mu z)^\gamma + B(1-\mu z)^{\gamma+\Delta},$$ where $0 < \Delta < 1.$ Such behaviour implies, at the coefficient level, 
$$s_n \sim \frac{A}{\Gamma(-\gamma)}\mu^nn^{-\gamma-1}+ \frac{B}{\Gamma(-\gamma-\Delta)}\mu^nn^{-\gamma-\Delta-1},$$ and for the ratios $$r_n \sim \mu\left ( 1-\frac{\gamma+1}{n}+\frac{const.}{n^{1+\Delta}}\right ).$$

From Fig. \ref{fig:r1} it is seen that $\mu \approx 12.4.$ Assuming this value, and estimating the gradient from the last plotted point, we find $g \approx -2.5.$ With $\mu = 12.3,$ we get $g=-2.1,$ and with $\mu = 12.5,$ we get $g=-2.9$ by this procedure, so it is clear that the estimate of $g$ is very sensitive to the estimate of $\mu.$

\begin{figure}[h!] 
   \centering
   \includegraphics[width=4in]{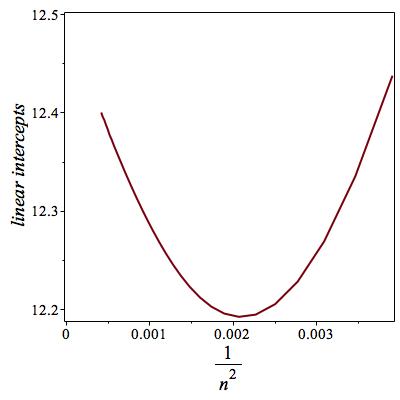} 
   \caption{Plot of intercepts of successive ratios against $1/n^2.$}
   \label{fig:int1}
\end{figure}

Calculating linear intercepts usually gives a more precise estimate of $\mu.$ One has $$n \cdot r_n -(n-1) \cdot r_{n-1} = \mu \left ( 1  + O \left ( \frac{1}{n^\Delta} \right ) \right ).$$ For a simple algebraic singularity, which means there is no confluent singularity at $z=z_c=1/\mu,$ so that $\Delta=1$ and the subdominant term is $O \left ( \frac{1}{n^2} \right ),$ and so convergence to $\mu$ is usually more rapid. However a plot of linear intercepts against $1/n^2,$ shown in Fig \ref{fig:int1} has gradient that changes sign for the last few  values of $n,$ and which may change sign again as $n$ increases, making it difficult to extrapolate, and strongly suggesting the presence of one or more confluent terms. It also implies that we would really need many more series coefficients in order to make more precise estimates of the critical parameters. It also reinforces the usefullness of the sequence extension procedure we have undertaken, as these approximate coefficients are essential to see this change of gradient. Despite these qualifications, a limiting value around $\mu=12.4,$ consistent with the value found by a simple ratio plot, seems plausible.

One can also calculate the gradient directly, from $$\frac{(r_n-r_{n-1})\cdot n(n-1)}{\mu} = g \left ( 1  + o \left ( \frac{1}{n} \right ) \right ).$$ Assuming the values $\mu=12.3,$ $\mu = 12.4,$ and $\mu=12.5$ we have plotted these estimators of $g$ against $1/n^2,$ as we don't know the correct sub-dominant exponent to use. Again one sees the necessity of estimating the last few terms, as otherwise the gradient change would not be observed, and a quite inaccurate estimate of $g$ would be obtained. As it is, we don't know how this plot will behave as $n$ increases, so cannot give any extrapolation with much confidence. However, if present trends continue, a value of $g \approx -2.5$ is plausible.

\begin{figure}[h] 
   \centering
   \includegraphics[width=4in]{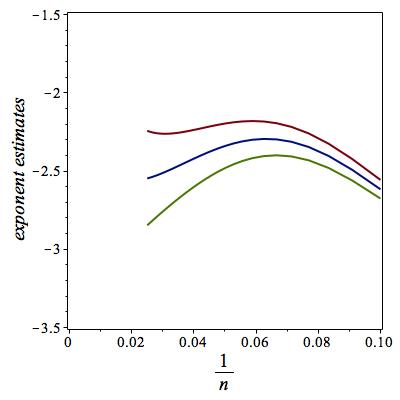} 
   \caption{Plot of estimators of exponent $g$ against $1/n^2.$ The top curve assumes $\mu=12.3,$ the middle curve assumes $\mu=12.4$ and the bottom curve assumes $\mu=12.5.$ }
   \label{fig:g}
\end{figure}

It is clear that the  $o \left ( \frac{1}{n} \right ) $ term is significantly affecting our extrapolation attempts. In an effort to address this, we make comparison with two similar problems whose asymptotics we have recently studied \cite{EPG15}. These are the corresponding problem of enumerating the number of permutations of various sizes that can be sorted by a deque, and by two stacks in parallel (2SIP). 

We have shown \cite{EPG15} that these two OGFs appear to have the same radius of convergence, which is quite accurately estimated as $x_c = 1/\mu_d \approx 0.120752497575574.$\footnote{This is more precise than the estimate in \cite{EPG15}, and is based on an analysis by Nathan Clisby of longer series.}  The generating function coefficients in the two cases are believed to behave as $d_n \sim const \cdot \mu_d^n \cdot n^{g_d},$ and  $p_n \sim const \cdot \mu_d^n \cdot n^{g_p},$
for deques and 2SIP respectively. Further, we have estimated that $g_d = -1.5$ and $g_p \approx -2.47327.$ All these data are based on an analysis of series of length 500 terms, so are vastly more reliable and precise than the estimates of the corresponding critical parameters in the current problem.

If we form the coefficient-by-coefficient quotients $s_n/d_n$ and $s_n/p_n,$ these will behave as $$s_d(n) = s_n/d_n \sim const \cdot \lambda^n \cdot n^{g_s-g_d}$$ and  $$s_p(n) = s_n/p_n \sim const \cdot \lambda^n \cdot n^{g_s-g_p},$$ respectively, where $\lambda = \mu/\mu_d.$ 

Now we can apply simple ratio analysis to the ratios $r_1(n) = s_d(n)/s_d(n-1),$ and $r_2(n) = s_p(n)/s_p(n-1).$ When plotted against $1/n,$ these should approach a common limit $\lambda,$ with gradient $\lambda(g_s-g_d)$ and $\lambda(g_s-g_p)$ respectively.

\begin{figure}[h] 
   \centering
   \includegraphics[width=4in]{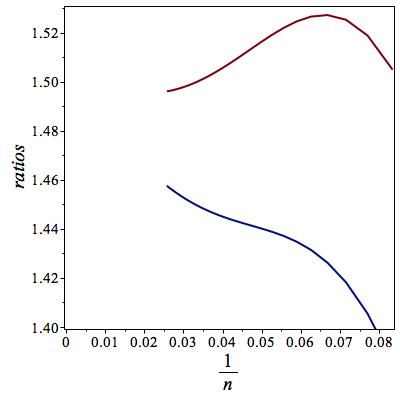} 
   \caption{Plot of coefficient ratios, as explained in text, against $1/n.$ Upper plot is for
two-stacks-in-series/two-stacks-in-parallel, lower plot is for two-stacks-in-series/deques. The common intersection point on the abscissa is estimated to be at 1.51}
   \label{fig:r2}
\end{figure}

\begin{figure}[h] 
   \centering
   \includegraphics[width=4in]{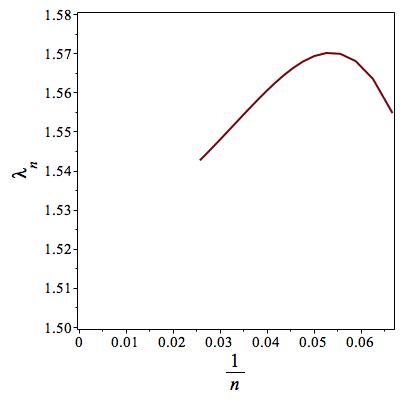} 
   \caption{Plot of $\lambda_n$ estimate, as explained in text, against $1/n.$ }
   \label{fig:inter}
\end{figure}

The ratio plots are shown in Fig \ref{fig:r2}, and it can be seen that the common limit is around $1.50-1.52,$ and that the gradients are of opposite sign. In fact the difference in gradients is $\lambda(g_p - g_d),$ and we know that $g_p - g_d \approx -0.97327.$ So we can tune the value of $\lambda$ to be consistent with this value,  as $$\lambda_n =\frac{n(r_1(n)-r_2(n))}{g_d-g_p} \left ( 1 + O\left( \frac{1}{n} \right ) \right ).$$ Plotting $\lambda_n$ against $1/n,$ shown in Fig. \ref{fig:inter}, we esimate $\lambda=1.51 \pm 0.01,$ which implies $\mu = 12.5 \pm 0.1,$ which is just consistent with previous analyses discussed above.

 This gives $\mu \approx 12.5,$ and $g_s \approx -2.5.$ Thus we take as our final estimates 
$\mu = 12.45 \pm 0.15$ and $g_s = -2.5 \pm 0.3,$ where the quoted errors are uncertainty estimates, and not in any sense rigorous error bounds.

Alternatively expressed, the OGF for two stacks in series behaves as $$ S(z) \sim const \cdot \left ( 1 - \mu \cdot z \right )^{-g_s-1} \approx  const \cdot \left ( 1 - 12.45  z \right )^{1.5}.$$

Our estimate of $\mu$ is of course consistent with the rigorous bounds given in \cite{AAL10}, which are $8.156 < \mu < 13.374.$ It is not inconceivable that the exponent could be the same as for two stacks in parallel, that is, 1.47327, but we have insufficient data to estimate the exponent with anything like this precision. This would correspond to the two problems being in the same universality class, when viewed from a statistical mechanical perspective.

Assuming the central estimates of both $\mu$ and the exponent $g$, one can estimate the amplitude by simple extrapolation. That is to say, if $s_n \sim a\cdot \mu^n \cdot n^g,$ then $a$ can be estimated by extrapolating the sequence $s_n/(\mu^n \cdot n^g)$ against $1/n.$ In this way we estimated $a \approx 0.008.$ Note however that this estimate is very sensitive to the estimates of both $\mu$ and $g.$ Writing the singular part of the generating function as $S(z) \sim A \cdot \left ( 1 - \mu \cdot z \right )^{-g-1},$ we have $A=a\Gamma(g+1) \approx 0.02.$

\section{Conclusion}
We have given an algorithm to generate the number of permutations of length $n$ sortable by two stacks in series. We have obtained the coefficients in the corresponding generating function up to and including permutations of length 19. We have used differential approximants to calculate the next 19 coefficients approximately, and the next 30 ratios of successive terms, and then analysed the extended series. In this way we have estimated the asymptotics of the generating function. We believe that the series length needs to be at least doubled in order to get much more significant accuracy in estimates of the critical parameters.

It is a source of some frustration that this problem appears to be so much harder than the corresponding problem of two stacks in parallel, for which an exact solution \cite{BMA14} is now available, as well as more than 1000 terms in the generating function. 

\section*{Acknowledgement}
We wish to thank Andrew Conway for many helpful discussions and an independent check of the first 10 coefficients of our expansion. We also thank Nathan Clisby for providing his variable precision version of the differential approximant program, which we used to extend the series, and Vince Vatter for corrections to references in an earlier version of this article.

\end{document}